\documentclass[twoside,reqno,A4]{amsart}

\usepackage{latexsym,rotate,eucal,cite,url}
\usepackage{amsmath,amsthm,amssymb,amsxtra} 

\newcommand{\puab}{\hspace*{-1.5mm}}

\newcommand{\dst}{\displaystyle}

\newcommand{\pav}[1]{\lfloor{#1}\rfloor}

\newcommand{\pavv}[1]{\left\lfloor{#1}\right\rfloor}

\newcommand{\ang}[1]{\langle{#1}\rangle}

\newcommand{\mb}[1]{\mathbb{#1}}

\newcommand{\hyp}[3]{\left[\puab\ba{#1}#2\\#3\ea\puab\right]} 
\newcommand{\hyq}[4]{\left[\puab\ba{#1}#3\\#4\ea{\puab\Big|#2}\right]}

\newcommand{\binm}{\binom}

\newcommand{\be}{\begin{equation}}
\newcommand{\ee}{\end{equation}}
\newcommand{\ba}{\begin{array}}
\newcommand{\ea}{\end{array}}
\newcommand{\bmn}{\begin{eqnarray}}
\newcommand{\emn}{\end{eqnarray}}
\newcommand{\bnm}{\begin{eqnarray*}}
\newcommand{\enm}{\end{eqnarray*}}
\newcommand{\bln}{\begin{subequations}}
\newcommand{\eln}{\end{subequations}}

\newcommand{\pq}[1]{\begin{equation}#1\end{equation}}
\newcommand{\pp}[2]{\begin{aligned}[#1]#2 
            \end{aligned}}  
 
\newcommand{\pmq}[1]{\begin{align}#1
            \end{align}}     
\newcommand{\pnq}[1]{\begin{align*}#1
            \end{align*}}    
\newcommand{\pnp}[2]{\begin{alignat*}{#1}#2
            \end{alignat*}}  
\newcommand{\centro}[1]
           {\begin{center}#1\end{center}}

\newcommand{\alp}{\alpha}
\newcommand{\bet}{\beta}
\newcommand{\gam}{\gamma}

\newcommand{\lam}{\lambda}

\newcommand{\Gam}{\Gamma}

\newtheorem{thm}{Theorem}
\newtheorem{lemm}[thm]{Lemma}
\newtheorem{corl}[thm]{Corollary}
\newtheorem{prop}[thm]{Proposition}

\newtheorem{entry}{Entry}

\newcommand{\referxy}[4]{\bibitem{kn:#1}{#2,}~\emph{#3,}~{#4.}}	
\newcommand{\cito}[1]{\cite{kn:#1}}	
\newcommand{\citu}[2]{\cite[#2]{kn:#1}}

\begin{document}
\title{Alternating Convolutions of Catalan Numbers}
\author{Wenchang Chu}
\address{School of Mathematics and Statistics\newline
         Zhoukou Normal University\newline
         Zhoukou (Henan), P.~R.~China\newline
\newline\emph{Corresponding address}:\newline
        Department of Mathematics and Physics\newline
     	University of Salento (P.~O.~Box~193) \newline
      	73100 Lecce, ~ Italy}
\email{chu.wenchang@unisalento.it}
\thanks{Email address: chu.wenchang@unisalento.it}
\subjclass[2010]{Primary 05A10, Secondary 33C15}
\keywords{Alternating convolution; Catalan number;
          Binomial coefficient; Hypergeometric series; 
          Product formula of confluent hypergeometric series}


\begin{abstract}
A new class of alternating convolutions concerning binomial 
coefficients and Catalan numbers are evaluated in closed forms.
\end{abstract}

\maketitle\thispagestyle{empty}
\markboth{Wenchang Chu}{Alternating Convolutions of Catalan Numbers}

\section{Introduction and Motivation}
The Catalan numbers 
\[C_n=\frac1{n+1}\binm{2n}{n}
\quad\text{for}\quad n\in\mb{N}_0\]
are probably the most frequently encountered sequence
in mathematics. There exist numerous interpretations
in enumerative combinatorics and remarkable identities
about them that can be found in the monographs by 
Koshy~\cito{koshy}, Roman~\cito{roman} and Stanley~\cito{stanley}
as well as in \cite{kn:chu2017b,kn:chu2018b,kn:chu2016f}.
For instance, these numbers satisfy the nonlinear recurrence relation
\[C_{n+1}=\sum_{k=0}^nC_kC_{n-k}\]
and the Touchard identity
\[C_{n+1}=\sum_{k=0}^{\pav{\frac{n}2}}2^{n-2k}\binm{n}{2k}C_k.\]
Here and forth, $\pav{x}$ denotes the integer part for the real number $x$.
For $m\in\mb{N}_0$ and $i,~j\in\mb{Z}$, we shall utilize the notation 
``$i\equiv_mj$" for ``$i$ is congruent to $j$ modulo $m$". The logical 
function $\chi$ is defined for brevity by $\chi(\text{true})=1$ and 
$\chi(\text{false})=0$. 

Recently, Miki\'c\citu{miki}{2019} found, by combinatorial bijections,
the following unusual convolution identities:
\pmq{
&\sum_{k=0}^n(-1)^k\binm{n}{k}
C_{k}C_{n-k}
=\frac{2\chi(n\equiv_20)}{n+2}
\binm{n}{\pav{\frac{n}2}}^2,\\
&\sum_{k=0}^n(-1)^k\binm{n}{k}
\binm{2n-2k}{n-k}C_{k}
=\binm{n}{\pav{\frac{n}2}}^2.}
Prodinger~\cito{prodg} provided different proofs 
by making use of Zeilberger's algorithm, Dixon's 
formula and its variants.
Let $(x)_n$ be the Pochhammer symbol given by
\[(x)_0=1\quad\text{and}\quad
(x)_n=x(x+1)\cdots(x+n-1)
\quad\text{for}\quad n\in\mb{N}.\]
The objective of this paper is to show the following generalizations.
\begin{thm}[$n,~\lam\in\mb{N}_0$]\label{thm=A}
\[\sum_{k=0}^n(-1)^k\binm{n}{k}
C_{k+\lam}C_{n-k+\lam}
=\frac{\lam!\chi(n\equiv_20)}{(2+n)_\lam}
\binm{2\lam}{\lam}\binm{n}{\pav{\frac{n}2}}
C_{\lam+\pav{\frac{n}2}}.\]
\end{thm}
\begin{thm}[$n,~\lam\in\mb{N}_0$]\label{thm=B}
\[\sum_{k=0}^n(-1)^k\binm{n}{k}
\binm{2n-2k+2\lam}{n-k+\lam}C_{k+\lam}
=\frac{\lam!}{(2+n)_\lam}
\binm{2\lam}{\lam}
\binm{n}{\pav{\frac{n}2}}
\binm{n+2\lam}{\lam+\pav{\frac{n}2}}.\]
\end{thm}

Three further binomial identities of alternating convolutions will also
be established.
\begin{thm}[$n,~\lam\in\mb{N}_0$]\label{thm=C}
\pnq{&\sum_{k=0}^n(-1)^k\binm{n}{k}
\binm{2k+2\lam}{k+\lam}\binm{2n-2k+2\lam}{n-k+\lam}\\
&~=\frac{\lam!\chi(n\equiv_20)}{(1+n)_\lam}
\binm{2\lam}{\lam}\binm{n}{\pav{\frac{n}2}}
\binm{2\lam+n}{\lam+\pav{\frac{n}2}}.}
\end{thm}
\begin{thm}[$n,~\lam\in\mb{N}_0$]\label{thm=D}
\[\pp{c}{&\sum_{k=0}^n(-1)^k\binm{n}{k}
\binm{2k+2\lam}{k+\lam}\binm{2n-2k+2\lam}{n-k+\lam}(n-k+\lam)\\
&~=\frac{\lam!}{(n)_{\lam}}
\binm{n}{\pav{\frac{n}2}}
\binm{2\lam}{\lam}\binm{2\lam+n}{\lam+\pav{\frac{n}2}}
\times\begin{cases}
\frac{n(2\lam+n)}{2(\lam+n)},&n\equiv_20;\\[2mm]
\frac{(n+1)(2\lam+n+1)}{2(\lam+n)},&n\equiv_21.
\end{cases}}\]
\end{thm}
\begin{thm}[$n,~\lam,~\mu\in\mb{N}_0$]\label{thm=E}
\[\sum_{k=0}^n(-1)^k\binm{n}{k}
\frac{\dst\binm{2k+2\lam}{k+\lam}\binm{2n-2k+2\mu}{n-k+\mu}}
{\dst\binm{k+2\lam}{\lam}\binm{n-k+2\mu}{\mu}}
=\frac{\dst\binm{n}{\frac{n}2}\binm{n+\lam+\mu}{\frac{n}2}}
{\dst\binm{\lam+\frac{n}2}{\lam}\binm{\mu+\frac{n}2}{\mu}}
\chi(n\equiv_20).\]
\end{thm}

The rest of the paper will be organized as follows. As a preliminary,
we shall prove, in the next section, a product formula for confluent 
hypergeometric $_1F_1$ and $_2F_2$-series, which may serve as a 
counterpart of the product formulae due to Bailey~\cito{bai28}. 
By extracting the coefficients of $x^n$ across these hypergeometric 
equations, we derive, in Section~3, three binomial formulae of 
alternating sums, which contain the five summation theorems just
displayed as particular cases. Finally, the paper will end with Section~4,
where two equivalent integral formulae are proposed as problems.

\section{Products of Hypergeometric Series} 
Recall that the $\Gam$-function~(see~\citu{rain}{\S8}
for example) is given by the beta integral
\[\Gam(x)=\int_{0}^{\infty}u^{x-1}
e^{-u}\mathrm{d}u
\quad\text{for}\quad
\Re(x)>0.\]
It satisfies Euler's reflection property
\[\Gam(x)\times\Gam(1-x)
=\frac{\pi}{\sin\pi x}\]
and Legendre's  duplicate formula
\[\Gam(2x)=\Gam(x)\Gam(x+\tfrac12)\frac{2^{2x-1}}{\sqrt{\pi}}.\]

For the sake of brevity, we shall utilize the following
multiparameter expression
\[\Gam\hyp{c}{\alp,\bet,\cdots,\gam}{A,B,\cdots,C}
=\frac{\Gam(\alp)\Gam(\bet)\cdots\Gam(\gam)}
      {\Gam(A)\Gam(B)\cdots\Gam(C)}.\]
Analogously, the quotient of the Pochhammer symbol will be abbreviated to 
\[\hyp{c}{\alp,\bet,\cdots,\gam}{A,B,\cdots,C}_n
=\frac{(\alp)_n(\bet)_n\cdots(\gam)_n}{(A)_n(B)_n\cdots(C)_n}.\]

According to Bailey~\citu{bai35}{\S2.1}, 
the classical hypergeometric series reads as
\[_{p}F_{q}\hyq{r}{z}
    {a_{1},a_{2},\cdots,a_{p}}
    {b_{1},b_{2},\cdots,b_{q}} 
=\sum_{k=0}^{\infty}
\frac{z^k}{k!}
\hyp{c}{a_1,a_2,\cdots,a_p}{b_1,b_2,\cdots,b_q}_k.\]
When $p\le q$, this series is said confluent. In 1928, 
Bailey~\cito{bai28} found, among others, two product 
formulae for confluent hypergeometric series
\pmq{{_1F_1}\hyq{c}{x}{a}{c}
\times\label{1F1-dixon}
{_1F_1}\hyq{c}{-x}{a}{c}
&={_2F_3}\hyq{c}{\frac{x^2}4}
{a,~c-a}{c,\frac{c}2,\frac{1+c}2},\\
{_1F_1}\hyq{c}{x}{a}{2a}
\times\label{1F1-watson}
{_1F_1}\hyq{c}{-x}{c}{2c}
&={_2F_3}\hyq{c}{\frac{x^2}4}
{\frac{a+c}2,\frac{a+c+1}2\\[-3.5mm]}
{a+c,a+\frac12,c+\frac12}.}
They resemble the following beautiful product
formula (cf.~Bailey~\citu{bai35}{\S10.1}) 
discovered by Clausen one century earlier
\[{_2F_1}^2\hyq{c}{x}
{a,~c}{a+c+\frac12}
={_3F_2}\hyq{c}{x}
{a+c,~2a,~2c}{a+c+\frac12,2a+2c}.\]

By inserting an extra linear factor $\lam+k$ in \eqref{1F1-dixon},
we find the extended formula.
\begin{lemm}[$a,~c,~\lam\in\mb{R}$]\label{1F1+linear} 
\pnq{&{_1F_1}\hyq{c}{x}{a}{c}
\times{_2F_2}\hyq{r}{-x}{1+\lam,a}{\lam,c}\\
&={_3F_4}\hyq{r}{\frac{x^2}4}
{1+\lam,~a,~c-a}{\lam,~c,\frac{c}2,\frac{1+c}2}
-\frac{ax}{c\lam}
{_2F_3}\hyq{c}{\frac{x^2}4}
{1+a,~c-a}{c,\frac{1+c}2,\frac{2+c}2}.}
\end{lemm}

In particular for $\lam=c-1$, we get 
a variant formula of \eqref{1F1-dixon}:
\pq{\label{1F1-linear}\pp{c}{
{_1F_1}\hyq{c}{x}{a}{c-1}
\times{_1F_1}\hyq{c}{-x}{a}{c}
={_2F_3}\hyq{c}{\frac{x^2}4}
{a,~c-a}{c-1,\frac{c}2,\frac{1+c}2}&\\
-\frac{ax}{c(c-1)}
{_2F_3}\hyq{c}{\frac{x^2}4}
{1+a,~c-a}{c,\frac{1+c}2,\frac{2+c}2}&.}}

\emph{Proof of Lemma~\ref{1F1+linear}}. \
By means of the linear relation
\[\frac{\lam+k}{k}=\frac{\lam-a}{\lam}+\frac{a+k}{\lam}\]
it is not hard to get the contiguous relation
\pnq{{_4F_3}\hyq{c}{1}{a,~c,~e,~1+\lam}{1+a-c,1+a-e,\lam}
&=\frac{\lam-a}{\lam}{_3F_2}\hyq{c}{1}{a,~c,~e}{1+a-c,1+a-e}\\
&\times\frac{a}{\lam}{_3F_2}\hyq{c}{1}{1+a,~c,~e}{1+a-c,1+a-e},}
where the condition $\Re(\frac{a}2-c-e)>0$ is provided for convergence.
Evaluating the former $_3F_2$-series by the Dixon formula 
(cf.~Bailey~\citu{bai35}{\S3.1})
\[{_3F_2}\hyq{c}{1}{a,~c,~e}{1+a-c,1+a-e}
=\Gam\hyp{c}
{1+a-c,1+a-e,1+\frac{a}2,1+\frac{a}2-b-c}
{1+\frac{a}2-c,1+\frac{a}2-e,1+a,1+a-c-e}\]
and the latter $_3F_2$-series by ``$D_{-1,-1}$" 
due to the author~\citu{chu12a}{Example~18}
\pnq{{_3F_2}\hyq{c}{1}{1+a,~c,~e}{1+a-c,1+a-e}
&=\frac{2^{2a-2c-2e-1}}{\pi}
\Gam\hyp{c}{1+a-c,1+a-e}{1+a-2c,1+a-2e}\\
&\pp{t}{\times\bigg\{\Gam&\hyp{c}{\tfrac{1+a}2,\tfrac{2+a}2-c,\tfrac{2+a}2-e,\tfrac{1+a}2-c-e}{1+a,~1+a-c-e}\\
+~\Gam&\hyp{c}{\tfrac{2+a}2,\tfrac{1+a}2-c,\tfrac{1+a}2-e,\tfrac{2+a}2-c-e}{1+a,~1+a-c-e}\bigg\}}}
and then simplifying the result, we get the expression
\pq{\label{4F3+linear}\pp{c}{
{_4F_3}\hyq{c}{1}{a,~c,~e,~1+\lam}{1+a-c,1+a-e,\lam}
&=\Gam\hyp{c}{1+a-c,1+a-e}{a,1+a-c-e}\\
\times\bigg\{\frac1{2\lam}
\Gam\hyp{c}{\frac{1+a}2,\frac{1+a}2-c-e}
{\frac{1+a}2-c,\frac{1+a}2-e}
&+\frac{2\lam-a}{4\lam}
\Gam\hyp{c}{\frac{a}2,\frac{2+a}2-c-e}
{\frac{2+a}2-c,\frac{2+a}2-e}\bigg\}.}}

In particular, when the series is terminated 
by $a=-n$ with $n\in\mb{N}_0$, we have
\pq{\label{4F3-linear}{_4F_3}\hyq{c}{1}
{-n,~c,~e,~1+\lam}{1-c-n,1-e-n,\lam}
=\hyp{c}{-n,1-c-e-n}{1-c-n,1-e-n}_{\pavv{\frac{n}2}}
\puab\times\begin{cases}
\tfrac{2\lam+n}{2\lam},&n\equiv_20;\\[2mm]
-\frac{1+n}{2\lam},&n\equiv_21.
\end{cases}}

Now we turn to examine, by letting $i+k=n$, the product
\pnq{{_1F_1}\hyq{c}{x}{a}{c}
\times{_2F_2}\hyq{r}{-x}{1+\lam,a}{\lam,c}
&=\sum_{i=0}^{\infty}
\frac{(a)_i}{i!(c)_i}
\sum_{k=0}^{\infty}(-1)^k
\frac{\lam+k}{\lam}
\frac{(a)_k}{k!(c)_k}x^{k+i}\\
&=\sum_{n=0}^{\infty}
\frac{(a)_nx^n}{n!(c)_n}
\sum_{k=0}^{n}
\hyp{rl}{-n,~a,&1-c-n,~\lam+1}
{1,~c,&1-a-n,~\lam}_k.}
Writing the last sum in terms of ${_4F_3}$-series
and then evaluating it by \eqref{4F3-linear} 
\[{_4F_3}\hyq{c}{1}
{-n,~a,~1-c-n,~\lam+1}
{~c,~1-a-n,~\lam}
=\hyp{c}{-n,c-a}{c,1-a-n}_{\pavv{\frac{n}2}}
\times\begin{cases}
\tfrac{2\lam+n}{2\lam},&n\equiv_20;\\[2mm]
-\frac{1+n}{2\lam},&n\equiv_21;
\end{cases}\]
we confirm, in view of the parity of $n$,
the product formula in Lemma~\ref{1F1+linear}.\qed

\section{Binomial Convolution Formulae} 
By extracting the coefficient of $x^n$ across the equations
\eqref{1F1-dixon}, \eqref{1F1-watson} and \eqref{1F1-linear}
on hypergeometric products, we find the following three identities.
\begin{prop}[$n\in\mb{N}_0$ and $a,~c\in\mb{R}$]\label{pp=A}
\[\sum_{k=0}^n(-1)^k\binm{n}{k}
\frac{(a)_k(a)_{n-k}}{(c)_k(c)_{n-k}}
=\frac{n!}{(c)_n}\hyp{c}{a,c-a}{1,~c}_{\pav{\frac{n}2}}
\chi(n\equiv_20).\]
\end{prop}

\begin{prop}[$n\in\mb{N}_0$ and $a,~c\in\mb{R}$]\label{pp=B}
\[\sum_{k=0}^n(-1)^k\binm{n}{k}
\frac{(a)_k(c)_{n-k}}{(2a)_k(2c)_{n-k}}
=\hyp{c}{1,a+c}{2a,~2c}_{n}
\hyp{c}{a,~c}{1,a+c}_{\pav{\frac{n}2}}
\chi(n\equiv_20).\]
\end{prop}

\begin{prop}[$n\in\mb{N}_0$ and $a,~c\in\mb{R}$]\label{pp=C}
\[\sum_{k=0}^n(-1)^k\binm{n}{k}
\frac{(a)_k(a)_{n-k}}{(c)_k(c-1)_{n-k}}
=\frac{n!}{(c-1)_{n+1}}\hyp{c}{a,c-a}{1,c}_{\pav{\frac{n}2}}
\times\begin{cases}
c+\dfrac{n-2}{2},&n\equiv_20;\\[2mm]
a+\dfrac{n-1}{2},&n\equiv_21.
\end{cases}\]
\end{prop}

Expressing the quotients of shifted factorials in terms 
of binomial coefficients
\pnp{4}{
\frac{(\frac12+\lam)_k}{(1+\lam)_k}
&=\frac{\binm{2k+2\lam}{k+\lam}}{4^k\binm{2\lam}{\lam}},
\quad&&\frac{(\frac32+\lam)_k}{(1+\lam)_k}
&&=\frac{\binm{2k+2\lam}{k+\lam}(1+2k+2\lam)}{4^k\binm{2\lam}{\lam}(1+2\lam)};\\
\frac{(\frac12+\lam)_k}{(2+\lam)_k}
&=\frac{C_{k+\lam}}{4^kC_{\lam}}, 
\quad&&\frac{(\frac32+\lam)_k}{(2+\lam)_k}
&&=\frac{\binm{1+2k+2\lam}{k+\lam}}{4^k\binm{1+2\lam}{\lam}};}
we can confirm the five identities anticipated in the introduction as follows:
\begin{itemize}
\item Theorem~\ref{thm=A}: $a=\frac12+\lam$ and $c=2+\lam$ in Proposition~\ref{pp=A}.
\item Theorem~\ref{thm=B}: $a=\frac12+\lam$ and $c=2+\lam$ in Proposition~\ref{pp=C}.
\item Theorem~\ref{thm=C}: $a=\frac12+\lam$ and $c=1+\lam$ in Proposition~\ref{pp=A}.
\item Theorem~\ref{thm=D}: $a=\frac12+\lam$ and $c=1+\mu$ in Proposition~\ref{pp=C}.
\item Theorem~\ref{thm=E}: $a=\frac12+\lam$ and $c=2+\mu$ in Proposition~\ref{pp=B}.
\end{itemize}

Furthermore, we can derive four ``reciprocal formulae" of those displayed 
in Theorems~\ref{thm=A}, \ref{thm=B}, \ref{thm=C} and~\ref{thm=D}.  
\begin{corl}[$a=2+\lam$ and $c=\frac12+\lam$ 
      in Proposition~\ref{pp=A}: $n,~\lam\in\mb{N}_0$]
\[\sum_{k=0}^n(-1)^k\binm{n}{k}
\frac{(n-1)(n-3)C^2_{\lam}}{C_{k+\lam}C_{n-k+\lam}}
=\frac{3C_{\lam}\binm{2\lam}{\lam}\binm{n}{\pav{\frac{n}2}}\chi(n\equiv_20)}
      {C_{\lam+\pav{\frac{n}2}}\binm{\lam+n}{\lam}\binm{2\lam+2n}{\lam+n}}.\]
\end{corl}
\begin{corl}[$a=2+\lam$ and $c=\frac32+\lam$ 
      in Proposition~\ref{pp=C}: $n,~\lam\in\mb{N}_0$]
\[\sum_{k=0}^n(-1)^k\binm{n}{k}
\frac{n\binm{1+2\lam}{\lam}C_{\lam}}
     {\binm{1+2k+2\lam}{k+\lam}C_{n-k+\lam}}
=\frac{(1+n+2\lam)C_{\lam}\binm{1+2\lam}{\lam}\binm{n}{\pav{\frac{n}2}}}
      {\binm{\lam+n+1}{n}\binm{2\lam+2n}{\lam+n}\binm{1+2\lam+n}{\lam+\pav{\frac{n}2}}}
\times\begin{cases}
\frac{n}{1-n},&n\equiv_20;\\[2mm]
\frac{1+n}{2-n},&n\equiv_21.
\end{cases}\]
\end{corl}
\begin{corl}[$a=1+\lam$ and $c=\frac12+\lam$ 
      in Proposition~\ref{pp=A}: $n,~\lam\in\mb{N}_0$]
\[\sum_{k=0}^n(-1)^k\binm{n}{k}
\frac{(1-n)\binm{2\lam}{\lam}^2}{\binm{2k+2\lam}{k+\lam}\binm{2n-2k+2\lam}{n-k+\lam}}
=\frac{\binm{2\lam}{\lam}^2\binm{n}{\pav{\frac{n}2}}\chi(n\equiv_20)}
      {\binm{\lam+n}{n}\binm{2\lam+2n}{\lam+n}\binm{2\lam+n}{\lam+\pav{\frac{n}2}}}.\]
\end{corl}
\begin{corl}[$a=1+\lam$ and $c=\frac32+\lam$ 
      in Proposition~\ref{pp=C}: $n,~\lam\in\mb{N}_0$]
\pnq{&\sum_{k=0}^n(-1)^k\binm{n}{k}
\frac{n(1+2\lam)(1+2n+2\lam)\binm{2\lam}{\lam}^2}
{(1+2k+2\lam)\binm{2k+2\lam}{k+\lam}\binm{2n-2k+2\lam}{n-k+\lam}}\\
&~=\frac{(1+2\lam)\binm{2\lam}{\lam}^2\binm{n}{\pav{\frac{n}2}}}
      {\binm{\lam+n}{n}\binm{2\lam+2n}{\lam+n}\binm{2\lam+n}{\lam+\pav{\frac{n}2}}}
\times\begin{cases}
n,&n\equiv_20;\\[2mm]
n+1,&n\equiv_21.
\end{cases}}
\end{corl}


\section{Integral representations} \
According to the expressions of the beta integrals (cf.~\cito{renzo})
\pnq{&\binm{2m}{m}=\frac{2^{2m}}{\pi}\bet(\tfrac12,\tfrac12+m)
=\frac{2^{2m}}{\pi}\int_0^1\frac{x^{m-\frac12}}{\sqrt{1-x}}\text{d}x,\\
&C_m=\frac{2^{1+2m}}{\pi}\bet(\tfrac32,\tfrac12+m)
=\frac{2^{1+2m}}{\pi}\int_0^1y^{m-\frac12}\sqrt{1-y}\text{d}y;}
we can reformulate the sum in Theorem~\ref{thm=A} as follows
\pnq{&\sum_{k=0}^n(-1)^k\binm{n}{k}
C_{k+\lam}C_{n-k+\lam}\\
&~=\frac{4^{1+n+2\lam}}{\pi^2}
\int_0^1\int_0^1
(xy)^{\lam-\frac12}
\sqrt{(1-x)(1-y)}
\sum_{k=0}^n(-1)^k\binm{n}{k}
x^{n-k}y^k\text{d}x\text{d}y\\
&~=\frac{4^{1+n+2\lam}}{\pi^2}
\int_0^1\int_0^1
(xy)^{\lam-\frac12}(x-y)^n
\sqrt{(1-x)(1-y)}\text{d}x\text{d}y.}
Then we get the following integral formula equivalent to Theorem~\ref{thm=A}.
\begin{corl}[$n,~\lam\in\mb{N}_0$]
\[\int_0^1\int_0^1
(xy)^{\lam-\frac12}(x-y)^n
\sqrt{(1-x)(1-y)}\text{d}x\text{d}y
=\frac{\pi^2~\lam!~\chi(n\equiv_20)}
      {4^{1+n+2\lam}(2+n)_\lam}
\binm{2\lam}{\lam}\binm{n}{\pav{\frac{n}2}}
C_{\lam+\pav{\frac{n}2}}.\]
\end{corl}

The sum in Theorem~\ref{thm=B} can analogously be manipulated: 
\pnq{&\sum_{k=0}^n(-1)^k\binm{n}{k}
\binm{2n-2k+2\lam}{n-k+\lam}C_{k+\lam}\\
&~=\frac{2^{1+2n+4\lam}}{\pi^2}
\int_0^1\int_0^1
(xy)^{\lam-\frac12}
\sqrt{\frac{1-y}{1-x}}
\sum_{k=0}^n(-1)^k\binm{n}{k}
x^{n-k}y^k\text{d}x\text{d}y\\
&~=\frac{2^{1+2n+4\lam}}{\pi^2}
\int_0^1\int_0^1
(xy)^{\lam-\frac12}(x-y)^n
\sqrt{\frac{1-y}{1-x}}
\text{d}x\text{d}y,}
which leads us to another integral formula.
\begin{corl}[$n,~\lam\in\mb{N}_0$]
\[\int_0^1\int_0^1
(xy)^{\lam-\frac12}(x-y)^n
\sqrt{\frac{1-y}{1-x}}
\text{d}x\text{d}y
=\frac{\pi^2~\lam!}{2^{1+2n+4\lam}(2+n)_\lam}
\binm{2\lam}{\lam}\binm{n}{\pav{\frac{n}2}}
\binm{n+2\lam}{\lam+\pav{\frac{n}2}}.\]
\end{corl}

Two further integral formulae corresponding to Theorems~\ref{thm=C}
and ~\ref{thm=D} can be produced in a similar manner. Finally,
an intriguing question is how to evaluate these integrals directly?



\end{document}